\documentclass[reqno,twoside,11pt,english]{amsart}
\usepackage{amsmath,amsfonts,amssymb,amsthm,epsfig}
\usepackage{amstext}
\usepackage{esint}
\makeatletter
\usepackage{color}
\usepackage[final,allcolors=blue,colorlinks=true]{hyperref}

\def\R{\mathbb R}
\def\N{\mathbb N}

\def\S{\mathbb S}

\def\al{\alpha}

\def\de{\delta}

\def\la{\lambda}
\def\si{\sigma}

\def\var{\varphi}

\def\na{\nabla}
\def\Om{\Omega}  
\def\De{\Delta}      

\def\cal{\mathcal}
\def\L{\mathcal L}                                       
\def\wq{\infty}
\def\pa{\partial}

\def\rad{\text{\rm rad}}

\def\Ker{\text{\rm Ker}}
\def\span{\text{\rm span}}
\newcommand{\D}{{\rm d}}

\newcommand{\medint}{-\kern -,375cm\int}         
\newcommand{\medintinrigo}{-\kern -,315cm\int}
\newcommand\esssup{\text{\rm \,esssup\,}}

\numberwithin{equation}{section}
\newtheorem{theorem}{Theorem}[section]

\newtheorem{corollary}[theorem]{Corollary}

\newtheorem{definition}[theorem]{Definition}

\newtheorem{proposition}[theorem]{Proposition}

\theoremstyle{definition}

\makeatother

\usepackage{babel}
\begin{document}
\title[Positive solutions of Kirchhoff equations]{Uniqueness and nondegeneracy of positive solutions to a class of  Kirchhoff  equations in $\R^3$}

 \author[G. Li, S. Peng and C.-L. Xiang]{Gongbao Li, Shuangjie Peng and Chang-Lin Xiang}

\address[Gongbao Li]{School of Mathematics and Statistics and Hubei Key Laboratory of Mathematical Sciences, Central China Normal University, Wuhan,  430079, P.R. China}
\email[]{ligb@mail.ccnu.edu.cn}
\address[Shuangjie Peng]{School of Mathematics and Statistics and Hubei Key Laboratory of Mathematical Sciences, Central China Normal University,  Wuhan,  430079, P.R. China}
\email[]{sjpeng@mail.ccnu.edu.cn}
\address[Chang-Lin Xiang]{School of Information and Mathematics, Yangtze University, Jingzhou 434023, P.R. China,  and
University of Jyvaskyla, Department of Mathematics and Statistics, P.O. Box 35, FI-40014 University of Jyvaskyla, Finland}
\email[]{xiang\_math@126.com}
\thanks{The first named author Li was supported by  Natural Science Foundation of China (Grant No. 11371159), Hubei Key Laboratory of Mathematical Sciences and Program for Changjiang Scholars and Innovative Research Team in University \# IRT13066.
 The corresponding author Xiang is financially supported by  the Academy of Finland, project 259224.}

\begin{abstract}
In this paper, we establish a type of  uniqueness and  nondegeneracy results for
 positive solutions to
  the following nonlocal Kirchhoff  equations
  \begin{eqnarray*}
-\left(a+b\int_{\R^{3}}|\na u|^{2}\D x\right)\De u+u=|u|^{p-1}u &
& \text{in }\R^{3},
\end{eqnarray*} where $a,b$ are positive constants and $1<p<5$. Before this paper, it seems that there have no this type of results even on positive ground states solutions to Kirchhoff type equations, much less on general positive solutions. To overcome the difficulty brought by the nonlocality,  some new observation on Kirchhoff equations is found, and some related theories on classical Schr\"odinger equations are applied.
\end{abstract}

\maketitle

{\small
\keywords {\noindent {\bf Keywords:} Kirchhoff equations; Nonlocality; Positive solutinos; Uniqueness; Nondegeneracy;}
\smallskip
\newline
\subjclass{\noindent {\bf 2010 Mathematics Subject
Classification:}  35A02, 35B20, 35J61}
}
\bigskip

\section{Introduction and main results}

\subsection{Introduction}

Let $a,b>0$ and $1<p<5$ be positive constants. In this paper we consider the following typical Kirchhoff type equations
\begin{eqnarray}
-\left(a+b\int_{\R^{3}}|\na u|^{2}\D x\right)\De u+u=|u|^{p-1}u &  & \text{in }\R^{3},\label{eq: Kirchhoff}
\end{eqnarray}
where  $u$ is a real-valued measurable
function, $\na u=(\pa_{x_{1}}u,\pa_{x_{2}}u,\pa_{x_{3}}u)$
 and $\De=\sum_{i=1}^{3}\pa_{x_{i}x_{i}}$ is the usual Laplacian
operator in $\R^{3}$.

Eq. (\ref{eq: Kirchhoff}) and its variants have been studied extensively
in the literature. The interest for studying Kirchhoff type equations
is twofold: first, the interest comes from the physical background
of Kirchhoff type equations. Indeed, to extend the classical D'Alembert's
wave equations for free vibration of elastic strings, Kirchhoff \cite{Kirchhoff-1883} proposed
the following time dependent wave equation
\[
\rho\frac{\pa^{2}u}{\pa^{2}t}-\left(\frac{P_{0}}{h}+\frac{E}{2L}\int_{0}^{L}\left|\frac{\pa u}{\pa x}\right|^{2}\D x\right)\frac{\pa^{2}u}{\pa x^{2}}=0
\]
for the first time. Some early classical studies of Kirchhoff equations
can be found in Bernstein \cite{Bernstein-1940} and Pohozaev \cite{Pohozaev-1975}.
Much attention was received after Lions \cite{Lions-1978} introducing
an abstract functional framework to this problem. More interesting
results in this respect can be found in e.g. \cite{DAncona-Spagnolo-1992,Arosio-Panizzi-1996,Cingolani-Lazzo-1997}
and the references therein. Second, the interest stems from the nonlocality
of Kirchhoff type equations from a mathematical point of view. For
instance, the consideration of the stationary analogue of Kirchhoff's
wave equation leads to the Dirichlet problem
\begin{equation}
\begin{cases}
-\left(a+b\displaystyle\int_{\Om}|\na u|^{2}\right)\De u=f(x,u) & \text{ in }\Om,\\
u=0 & \text{on }\pa\Om,
\end{cases}\label{eq: Kirchhoff eq. on bdd domain}
\end{equation}
 where $\Om\subset\R^{3}$ is a bounded domain, and to equations of
type
\begin{eqnarray}
-\left(a+b\int_{\R^{3}}|\na u|^{2}\right)\De u=f(x,u) &  & \text{in }\R^{3},\label{eq: general Kirchhoff eq.}
\end{eqnarray}
respectively. In above two equations, $f$ denotes some nonlinear
functions, a typical example of which is given as in Eq. (\ref{eq: Kirchhoff}).
In both equations (\ref{eq: Kirchhoff eq. on bdd domain}) and (\ref{eq: general Kirchhoff eq.}),
note that the term $\left(\int|\na u|^{2}\D x\right)\De u$ depends
not only on the pointwise value of $\De u$, but also on the integral
of $|\na u|^{2}$ over the  domain. In this sense, equations (\ref{eq: Kirchhoff}),
(\ref{eq: Kirchhoff eq. on bdd domain}) and (\ref{eq: general Kirchhoff eq.})
are no longer the usual pointwise equality. This new feature brings
new mathematical difficulties that make the study of Kirchhoff type
equations particularly interesting. We refer to e.g. \cite{He-Zou-2009,He-Zou-2010,Perera-Zhang-2006,Shuai-2015,Zhang-Perera-2006}
and to e.g. \cite{Deng-Peng-Shuai-2015,Figueiredo et al-2014,Guo-2015,He-2016,He-Li-2015,He-Li-Peng-2014,Li-Li-Shi-2012,Li-Ye-2014,Lu-2015-arXiv,Ye-2015}
for mathematical researches on Kirchhoff type equations on bounded
domains and in the whole space, respectively. Nonlocal problems also appear in other  mathematical research fields. We refer the interested readers
to e.g. \cite{Fall-Valdinoci-2014,Frank-Lenzmann-2013,Frank-Lenzman-Silvestre-2016,Lenzmann-2009}  and to \cite{Lenzmann-2009,Xiang-CVPDE}
 for mathematical researches on fractional
type nonlocal Schr\"odinger equations and convolution type nonlocal Choquard equations, repectively.

In this paper, we are concerned about positive solutions of Eq. (\ref{eq: Kirchhoff}).
By a solution, we mean a function $u$ in $H^{1}(\R^{3})$ such that
\begin{eqnarray*}
\int_{\R^{3}}\left(\left(a+b\int_{\R^{3}}|\na u|^{2}\D x\right)\na u\cdot\na\var+u\var-|u|^{p-1}u\var\right)\D x=0, &  & \forall\,\var\in C_{0}^{\wq}(\R^{3}).
\end{eqnarray*}
 It is well known that (\ref{eq: Kirchhoff}) is the Euler-Lagrange
equation of the energy functional $I:H^{1}(\R^{3})\to\R$ defined
as
\[
I(u)=\frac{1}{2}\int_{\R^{3}}\left(a|\na u|^{2}+u^{2}\right)\D x+\frac{b}{4}\left(\int_{\R^{3}}|\na u|^{2}\D x\right)^{2}-\frac{1}{p+1}\int_{\R^{3}}|u|^{p+1}\D x
\]
for $u\in H^{1}(\R^{3})$. Thus critical point theories have been
devoted to find solutions for Eq. (\ref{eq: Kirchhoff}) and its
variants, see e.g.
\cite{He-Zou-2009,He-Zou-2010,He-Zou-2012,Li-Ye-2014,Ye-2015} and
the references therein. In particular, the existence of positive
solutions of Eq. (\ref{eq: Kirchhoff}) was obtained by looking for
 the so called ground states, which is defined as follows: Consider
the set of solutions to Eq. (\ref{eq: Kirchhoff}) and denote
\begin{equation}
m=\inf\left\{ I(v):v\in H^{1}(\R^{3})\text{ is a nontrivial solution to Eq. }(\ref{eq: Kirchhoff})\right\} .\label{eq: least energy}
\end{equation}
A nontrivial solution $u$ to Eq. (\ref{eq: Kirchhoff}) is called
a \emph{ground state} if
\[
I(u)=m.
\]
Since we focus on positive solutions of Eq. (\ref{eq: Kirchhoff}),
it is convenient to summarize some known results in the literature
before we proceed further.

\begin{proposition} \label{prop: Existence}Let $a,b>0$ be positive constants
and $1<p<5$. Let $m$ be the ground state energy defined as in
(\ref{eq: least energy}). Then, there exists a ground state of
(\ref{eq: Kirchhoff}) which is positive, and there holds
\[
m>0.
\]
Moreover, for any positive solution $u$, there hold

(1) (smoothness) $u\in C^{\wq}(\R^{3})$;

(2) (symmetry) there exists a decreasing function $v:[0,\wq)\to(0,\wq)$
such that $u=v(|\cdot-x_{0}|)$ for a point $x_{0}\in\R^{3}$;

(3) (Asymptotics) For any multiindex $\al\in\N^{n}$, there exist
constants $\de_{\al}>0$ and $C_{\al}>0$ such that
\begin{eqnarray*}
|D^{\al}u(x)|\le C_{\al}e^{-\de_{\al}|x|} &  & \text{for all }x\in\R^{3}.
\end{eqnarray*}
\end{proposition}

The existence of ground states of Eq. (\ref{eq: Kirchhoff}) is implied
by Proposition \ref{prop: Existence} of Ye \cite{Ye-2015}%
\footnote{This reference was brought to us by Ye.%
}, where more general existence results on Kirchhoff type
equations in $\R^{3}$ are obtained. In the special cases when
$3<p<5$ and $2<p<3$, the existence has also been proved by He and
Zou \cite{He-Zou-2012} and Li and Ye \cite{Li-Ye-2014},
respectively. In particular, in the papers of Ye \cite{Ye-2015} and  Li and Ye \cite{Li-Ye-2014}, to apply the Mountain Pass Lemma to find a ground state solution,  quite complicated manifolds were constructed in order to find a bounded Palais-Smale sequence. The fact that $m>0$ follows from Li and Ye
\cite[Lemma 2.8]{Li-Ye-2014}, see also Ye \cite{Ye-2015}. Other
properties follow easily from the theory of classical
Schr\"odinger equations (This will become clearer in view of the
Theorem \ref{thm: Uniqueness} below). For applications of Proposition
\ref{prop: Existence}, see e.g. He and Zou \cite{He-Zou-2012}, Li
and Ye \cite{Li-Ye-2014} and Ye \cite{Ye-2015} and the references
therein.

\subsection{Motivations and main results}

Proposition \ref{prop: Existence} provides a good understanding on ground
states of Eq. (\ref{eq: Kirchhoff}). However, we are still left an open
problem of uniqueness and nondegeneracy of the ground state, which turns out to be important to
 know further quantitative properties
of ground states when one studies more general and difficult
problems concerning Kirchhoff type equations, such as singular
perturbation problems related to Eq. (\ref{eq: Kirchhoff}).
Concerning uniqueness and nondegeneracy of ground states, there
exist several interesting results.  For instance, it is well known
\cite{Berestycki-Lions-1983-1,Chang et al-2007,Kwong-1989} that
the classical Schr\"odinger equation
\begin{eqnarray}\label{eq: classical Schrodinger equation}
-\De w+w=w^{q}, & w>0 & \text{in }\R^{N}
\end{eqnarray}
admits a unique positive solution (up to translations) which is
also nondegenerate (see the Definition \ref{def: Nondegeneracy}
below). The same results also hold for positive solutions to the
quasilinear Schr\"odinger equation
\begin{eqnarray}\label{eq: quasilinear Schrodinger equation}
-\Delta u-u\Delta|u|^{2}+u-|u|^{q-1}u=0 &  & \text{in }\R^{N},
\end{eqnarray}
 see e.g. \cite{Selvitella-2015,Shinji-Masataka-Tatsuya-2016,Xiang-2016},
and for ground states of the fractional Schr\"odinger equations ($0<s<1\le N$)
\begin{eqnarray*}
\left(-\De\right)^{s}w+w=w^{q}, & w>0 & \text{in }\R^{N},
\end{eqnarray*}
see e.g.
\cite{Fall-Valdinoci-2014,Frank-Lenzmann-2013,Frank-Lenzman-Silvestre-2016}.
In above three examples, $q$ is an index standing for the
nonlinearity of subcritical growth.   For a systematical research on
applications of nondegeneracy of ground states to perturbation
problems, we refer to Ambrosetti and Malchiodi
\cite{Ambrosetti-Malchiodi-Book} and the references therein. It is
also known that the uniqueness and nondegeneracy of ground states
are of fundamental importance when one deals with orbital
stability or instability of ground states. It mainly removes the
possibility that directions of instability come from the kernel of
$\L_{+}$ (see Definition \ref{def: Nondegeneracy}).
 The uniqueness and nondegeneracy
of ground states also play an important role in blow-up analysis for the corresponding standing wave solutions in the
corresponding time-dependent equations, see e.g. \cite{Frank-Lenzmann-2013,Frank-Lenzman-Silvestre-2016}
and the references therein.

Return to Kirchhoff equations. So far, there seems to have no result
on the uniqueness and nondegeneracy of ground states to equations
such as (\ref{eq: Kirchhoff}), much less on the general positive
solutions. Motivated by the fundamental importance of the uniqueness
and nondegeneracy of positive solutions and their numerous potential
applications as mentioned above, in this paper we aim to establish
the uniqueness and nondegeneracy of positive solutions to Eq. (\ref{eq: Kirchhoff}).
Our first main result reads as follows.

\begin{theorem} \label{thm: Uniqueness} Let $a,b>0$ be positive
constants and $1<p<5$. Then, positive solutions of Eq. (\ref{eq: Kirchhoff})
are unique up to translations. \end{theorem}

In particular, combining the symmetry result of Proposition \ref{prop:
Existence}, we infer that there exists a unique smooth positive
radial solution to Eq. (\ref{eq: Kirchhoff}) which decays
exponentially at infinity. Since ground states solutions of Eq.
(\ref{eq: Kirchhoff}) are of constant sign, we also infer from Theorem
\ref{thm: Uniqueness} that every positive solution  is exactly a
 positive ground state to Eq. (\ref{eq: Kirchhoff}).

We  remark that our proof of Theorem \ref{thm: Uniqueness} can also be seen as   a new proof for the existence of positive solutions to Eq. (\ref{eq: Kirchhoff}).  Indeed, by the proof of Theorem \ref{thm: Uniqueness}, we obtain an explicit expression for the positive solutions to Eq. (\ref{eq: Kirchhoff}).  Recall that to find  a positive ground state solution for Eq. (\ref{eq: Kirchhoff}),  Ye \cite{Ye-2015} and  Li and Ye \cite{Li-Ye-2014} constructed quite complicated manifolds so as to use the Mountain Pass Lemma. While in our proof, we
  derive all the positive solutions from a completely different way which is far more  elementary than that of
 He and Zou \cite{He-Zou-2010},  Ye \cite{Ye-2015} and  Li and Ye \cite{Li-Ye-2014}. Also, contrary to the different techniques applied to different ranges for the power $p$ in He and Zou \cite{He-Zou-2010}  and  Li and Ye \cite{Li-Ye-2014}, our approach is  unified for $p$ in the whole range $1<p<5$.

Our next main result concerns about nondegeneracy of the positive
solutions to Eq. (\ref{eq: Kirchhoff}) defined as follows.

\begin{definition}\label{def: Nondegeneracy} Let $u$ be a positive
solution of Eq. (\ref{eq: Kirchhoff}). We say that $u$
is nondegenerate in $H^{1}(\R^{3})$, if the following holds:\textbf{
}
\begin{equation}
\Ker\L_{+}=\span\left\{ \pa_{x_{1}}u,\pa_{x_{2}}u,\pa_{x_{3}}u\right\} ,\label{eq: Kernel}
\end{equation}
where $\L_{+}:L^{2}(\R^{3})\to L^{2}(\R^{3})$ is the linearized operator
around $u$ defined as
\[
\L_{+}\var=-\left(a+b\int_{\R^{3}}|\na u|^{2}\D x\right)\De\var+\var-pu^{p-1}\var-2b\left(\int_{\R^{3}}\na u\cdot\na\var\D x\right)\De u
\]
for $\var\in L^{2}(\R^{3})$. \end{definition}

It is easy to verify that $\L_{+}$ is a self-adjoint operator
acting on $L^{2}(\R^{3})$ with form domain $H^{1}(\R^{3})$. Note
that $\L_{+}$ is nonlocal due to the last term. In view of this
nonlocality, the following nondegenracy result is not obvious at all.

\begin{theorem} \label{thm: Nondegeneracy}Let $a,b>0$ be positive
constants and $1<p<5$. Then the positive solutions of Eq. (\ref{eq: Kirchhoff})
are nondegenerate in $H^{1}(\R^{3})$ in the sense of Definition \ref{def: Nondegeneracy}.
\end{theorem}

Theorem \ref{thm: Nondegeneracy} will be proved in section \ref{sec: nondegeneracy}
following the line of Ambrosetti and Malchiodi \cite[Chapter 4]{Ambrosetti-Malchiodi-Book}.

Some remarks are in order.  First, we
give a simple observation. It is easy to show that $\L_{+}$ is a
Fredholm operator of index zero, since the positive solution $u$
and its derivatives decay exponentially at infinity. Hence, if we
denote by ${\cal Z}$ the critical manifold consisting of all the
constant signed solutions of Eq. (\ref{eq: Kirchhoff}), then
${\cal Z}$ is nondegenerate in the sense of Ambrosetti and
Malchiodi \cite[Chapter 2]{Ambrosetti-Malchiodi-Book}. Second, we
note that our arguments are applicable to more general Kirchhoff
equations in $\R^{3}$ under suitable assumptions. However, in the
present paper we still restrict the research on Eq. (\ref{eq:
Kirchhoff}) due to its typicality.

Before closing this section, let us briefly show our idea and sketch the proofs. Recall that to deduce the uniqueness and nondegeneracy for positive solutions to the local Schr\"odinger equations (\ref{eq: classical Schrodinger equation}) and (\ref{eq: quasilinear Schrodinger equation}), corresponding ordinary differential equations are used. That is, to consider the ordinary differential equations
\begin{eqnarray*}
-\left(u_{rr}+\frac {N-1}{r}u_r\right)+u(r)-u^{p}(r)=0,&&r>0,
\end{eqnarray*}
and
\begin{eqnarray*}
-\left(u_{rr}+\frac {N-1}{r}u_r\right)-u(r)\left((u^2)_{rr}+\frac {N-1}{r}(u^2)_r\right)+u(r)-u^{p}(r)=0,&&r>0,
\end{eqnarray*} respectively, where $u_r$ is the derivative of $u$ with respect to $r$, see e.g. Kwong \cite{Kwong-1989} and Shinji et al. \cite{Shinji-Masataka-Tatsuya-2016}.
 Therefore,  to prove Theorem \ref{thm: Uniqueness}, it is quite natural to consider the corresponding ordinary differential equation to Eq. (\ref{eq: Kirchhoff})
\[-\left(a+b\int_0^{\wq}u_r^{ 2}(r)\D r\right)\left(u_{rr}+\frac {2}{r}u_r\right)+u(r)-u^{p}(r)=0\]
for $0<r<\wq$.  However, a further research shows that this idea is not so applicable due to the nonlocality of the term $\int_0^{\wq}u_r^{ 2}(r)\D r$. To overcome this difficulty, our key observation is  that the quantity $\int_0^{\wq}u_r^{ 2}(r)\D r$ is, in fact, independent of the positive solution $u$. Hence we conclude that the coefficient $ a+b\int_0^{\wq}u_r^{ 2}(r)\D r$ is no more than a positive constant that is independent of the given solution $u$. At this moment, we are allowed to apply the uniqueness result  of Kwong \cite{Kwong-1989}  on positive solutions to Eq. (\ref{eq: classical Schrodinger equation}) to prove Theorem \ref{thm: Uniqueness}.

 Next, to prove Theorem  \ref{thm: Nondegeneracy}, we apply the spherical harmonics to turn the problem into a system of ordinary differential equations.   It turns out that the key is to show that the problem $\L_+ \var=0$ has only a trivial radial solution. In other words, the key step is to show that the positive solution $u$ of Eq. (\ref{eq: Kirchhoff}) is nondegenerate in the subspace of radial functions of  $H^1(\R^3)$. To this end, again the above observation  plays an essential role. To be precise, write $c=a+b\int_0^{\wq}u_r^{ 2}(r)\D r$ and keep in mind that $c$ is a constant independent of $u$. Introduce an auxiliary operator $\cal{A}_u$  associated to $u$ by defining \[ {\cal A}_{u}\var=-c\De\var+\var-pu^{p-1}\var. \] Then solving the problem $\L_+ \var=0$, where $\var$ is radial,  is equivalent to solve \[
{\cal A}_{u}\var=2b\left(\int_{\R^{3}}\na u\cdot\na\var\D x\right)\De u.
\]
Since ${\cal A}_{u}$ is  the linearized operator of positive solutions to  Eq. (\ref{eq: classical Schrodinger equation}) up to a constant, the theory of the nondegeneracy of positive  solutions to Eq. (\ref{eq: classical Schrodinger equation}) are applicable, see Proposition \ref{thm: Kwong and Chang} and Proposition \ref{prop: property of Au}. Finishing this step, the remaining proof is standard. We refer the readers  to  the proof of Theorem \ref{thm: Nondegeneracy} for details.

Our notations are standard. $\N=\{0,1,2,\cdots\}$ denotes the set
of nonnegative integers. For any $1\le s\le\infty$, $L^{s}(\R^{3})$
is the Banach space of real-valued Lebesgue measurable functions $u$
such that the norm
\[
\|u\|_{s}=\begin{cases}
\left(\int_{\R^{3}}|u|^{s}\D x\right)^{1/s} & \text{if }1\le s<\infty\\
\esssup_{\R^{3}}|u| & \text{if }s=\infty
\end{cases}
\]
is finite. A function $u$ belongs to the Sobolev space $H^{1}(\R^{3})$
 if $u\in L^{2}(\R^{3})$ and its first order weak partial derivatives
also belong to $L^{2}(\R^{3})$. We equip $H^{1}(\R^{3})$
with the norm
\[
\|u\|_{H^{1}}=\sum_{\al\in\N^{3},|\al|\le 1}\|\pa^{\al}u\|_{2}.
\]
We also denote by $H_{\rad}^{1}(\R^{3})$ the subspace of radial Sobolev
functions in $H^{1}(\R^{3})$. For the properties of the Sobolev functions,
we refer to the monograph \cite{Ziemer}. By the usual abuse of notations,
we write $u(x)=u(r)$ with $r=|x|$ whenever $u$ is a radial function
in $\R^{3}$.

\section{Uniqueness of positive solutions\label{sec: Uniqueness}}

In this section we prove Theorem \ref{thm: Uniqueness}. Throughout
the following two sections, we denote by $Q\in H^{1}(\R^{3})$ the
unique positive radial function that satisfies
\begin{eqnarray}
-\De Q+Q=Q^{p} &  & \text{in }\R^{3}.\label{eq: Kwong}
\end{eqnarray}
 We refer to e.g. Berestycki and Lions \cite{Berestycki-Lions-1983-1}
and Kwong \cite{Kwong-1989} for the existence and uniqueness of $Q$,
respectively.

\begin{proof}[Proof of Theorem \ref{thm: Uniqueness}]Let $u\in H^{1}(\R^{3})$
be an arbitrary positive solution to Eq. (\ref{eq: Kirchhoff}). Write
$c=a+b\int_{\R^{3}}|\na u|^{2}\D x$ so that $u$ satisfies
\begin{eqnarray*}
-c\De u+u=u^{p} &  & \text{in }\R^{3}.
\end{eqnarray*}
Then, it is direct to verify that $u(\sqrt{c}(\cdot-t))$ solves Eq.
(\ref{eq: Kwong}) for any $t\in\R^{3}$. Thus, the uniqueness of
$Q$ implies that
\begin{eqnarray*}
u(x)=Q\left(\frac{x-t}{\sqrt{c}}\right), &  & x\in\R^{3},
\end{eqnarray*}
for some $t\in\R^{3}$. In particular, we obtain $\int_{\R^{3}}|\na u|^{2}\D x=\sqrt{c}\int_{\R^{3}}|\na Q|^{2}\D x$.
Substituting this equality into the definition of $c$ yields
\[
c=a+b\|\na Q\|_{2}^{2}\sqrt{c}.
\]
Since $c>0$, this equation is uniquely solved by
\begin{equation}
\sqrt{c}=\frac{1}{2}\left(b\|\na Q\|_{2}^{2}+\sqrt{b^{2}\|\na Q\|_{2}^{4}+4a}\right).\label{eq: value of c}
\end{equation}
As a consequence, we deduce that
\[
u(x)=Q\left(\frac{2(x-t)}{b\|\na Q\|_{2}^{2}+\sqrt{b^{2}\|\na Q\|_{2}^{4}+4a}}\right)
\]
for some $t\in\R^{3}$. At this moment, we can easily conclude that
the set
\[
{\cal M}=\left\{ Q\left(\frac{2(x-t)}{b\|\na Q\|_{2}^{2}+\sqrt{b^{2}\|\na Q\|_{2}^{4}+4a}}\right):t\in\R^{3}\right\}
\]
consists of all the positive solutions of Eq. (\ref{eq: Kirchhoff}).
The proof of Theorem \ref{thm: Uniqueness} is complete.\end{proof}

Note that (\ref{eq: value of c}) implies that the value of $c$ is independent
of the choice of positive solutions.

Before we end this section, we give a simple application of our uniqueness
result. Recall that $m$ is defined in (\ref{eq: least energy}) as
the ground state energy of the functional $I$. It is now available
to give an explicit expression of $m$ in terms of $a$, $b$ and
$\|\na Q\|_{2}$. We leave this to the interested readers since it
has no importance in the present paper. We point out that the following
result can be derived naturally.

\begin{corollary} The ground state energy $m$ is an isolated critical
value of $I$. \end{corollary}

\section{Nondegeneracy of positive solutions\label{sec: nondegeneracy}}

In this section we prove Theorem \ref{thm: Nondegeneracy}. We need
the following result.

\begin{proposition} \label{thm: Kwong and Chang} Let $1<p<5$ and
let $Q\in H^{1}(\R^{3})$ be the unique positive radial ground state
of Eq. (\ref{eq: Kwong}). Define the operator ${\cal A}:L^{2}(\R^{3})\to L^{2}(\R^{3})$
as
\[
{\cal A}\var=-\De\var+\var-pQ^{p-1}\var
\]
for $\var\in L^{2}(\R^{3})$. Then the following hold:

(1) $Q$ is nondegenerate in $H^{1}(\R^{3})$, that is,
\[
\Ker{\cal A}=\span\left\{ \pa_{x_{1}}Q,\pa_{x_{2}}Q,\pa_{x_{3}}Q\right\} ;
\]

(2) The restriction of ${\cal A}$ on $L_{\rad}^{2}(\R^{3})$ is one-to-one
and thus it has an inverse ${\cal A}^{-1}:L_{\rad}^{2}(\R^{3})\to L_{\rad}^{2}(\R^{3})$;

(3) ${\cal A}Q=-(p-1)Q^{p}$ and
\[
{\cal A}R=-2Q,
\]
where $R=\frac{2}{p-1}Q+x\cdot\na Q$.\end{proposition}

For a brief proof of (1), we refer to Chang et al. \cite[Lemma 2.1]{Chang et al-2007}
(see also the references therein); (2) is an easy consequence of (1)
since $Q$ is radial and ${\rm Ker}{\cal A}\cap L_{\rad}^{2}(\R^{3})=\emptyset$;
the last result can be obtained by a direct computation, see also
Eq. (2.1) of Chang et al. \cite{Chang et al-2007}.

Next, we introduce an auxiliary operator. Let $u$ be a positive solution
of Eq. (\ref{eq: Kirchhoff}). Since Eq. (\ref{eq: Kirchhoff}) is
translation invariant, we assume with no loss of generality that $u$
is radially symmetric with respect to the origin. Write $c=a+b\int_{\R^{3}}|\na u|^{2}\D x$.
Keep in mind that $c$ is a constant that is independent of the
choice of $u$ by (\ref{eq: value of c}). Then $u$ satisfies
\begin{eqnarray}
-c\De u+u-u^{p}=0 &  & \text{in }\R^{3}.\label{eq: auxiliary equation}
\end{eqnarray}
Define the auxiliary operator ${\cal A}_{u}:L^{2}(\R^{3})\to L^{2}(\R^{3})$
as
\[
{\cal A}_{u}\var=-c\De\var+\var-pu^{p-1}\var
\]
for $\var\in L^{2}(\R^{3})$. The following result on ${\cal A}_{u}$
follows easily from Proposition \ref{thm: Kwong and Chang}.

\begin{proposition} \label{prop: property of Au} ${\cal A}_{u}$
satisfies the following properties:

(1) The kernel of ${\cal A}_{u}$ is given by
\[
\Ker{\cal A}_{u}=\span\left\{ \pa_{x_{1}}u,\pa_{x_{2}}u,\pa_{x_{3}}u\right\} ;
\]

(2) The restriction of ${\cal A}_{u}$ on $L_{\rad}^{2}(\R^{3})$
is one-to-one and thus it has an inverse ${\cal A}_{u}^{-1}:L_{\rad}^{2}(\R^{3})\to L_{\rad}^{2}(\R^{3})$;

(3) ${\cal A}_{u}u=-(p-1)u^{p}$ and
\[
{\cal A}_{u}S=-2u,
\]
where $S=\frac{2}{p-1}u+x\cdot\na u$. \end{proposition}
\begin{proof}
Apply Proposition \ref{thm: Kwong and Chang} to $\tilde{u}$ defined
by $\tilde{u}(x)=u(\sqrt{c}x)=Q(x)$. We leave the details to the
interested readers.
\end{proof}
We will also use the standard spherical harmonics to decompose functions
in $H^{j}(\R^{N})$ for $j=0,1$, where $N=3$ (see e.g. Ambrosetti
and Malchiodi \cite[Chapter 4]{Ambrosetti-Malchiodi-Book}). So let
us introduce some necessary notations for the decomposition. Denote
by $\De_{\S^{N-1}}$ the Laplacian-Beltrami operator on the unit
$N-1$ dimensional sphere $\S^{N-1}$ in $\R^{N}$. Write
\begin{eqnarray*}
M_{k}=\frac{(N+k-1)!}{(N-1)!k!}\quad\forall\, k\ge0, & \text{and } & M_{k}=0\quad\forall\, k<0.
\end{eqnarray*}
Denote by $Y_{k,l}$, $k=0,1,\ldots$ and $1\le l\le M_{k}-M_{k-2}$,
the spherical harmonics such that
\[
-\De_{\S^{N-1}}Y_{k,l}=\la_{k}Y_{k,l}
\]
for all $k=0,1,\ldots$ and $1\le l\le M_{k}-M_{k-2}$, where
\begin{eqnarray*}
\la_{k}=k(N+k-2) &  & \forall\, k\ge0
\end{eqnarray*}
is an eigenvalue of $-\De_{\S^{N-1}}$ with
multiplicity$M_{k}-M_{k-2}$ for all $k\in\N$. In particular,
$\la_{0}=0$ is of multiplicity 1 with $Y_{0,1}=1$, and
$\la_{1}=N-1$ is of multiplicity $N$ with $Y_{1,l}=x_{l}/|x|$ for
$1\le l\le N$.

Then for any function $v\in H^{j}(\R^{N})$, we have the decomposition
\[
v(x)=v(r\Om)=\sum_{k=0}^{\wq}\sum_{l=1}^{M_{k}-M_{k-2}}v_{kl}(r)Y_{kl}(\Om)
\]
with $r=|x|$ and $\Om=x/|x|$, where
\begin{eqnarray*}
v_{kl}(r)=\int_{\S^{N-1}}v(r\Om)Y_{kl}(\Om)\D\Om &  & \forall\, k,l\ge0.
\end{eqnarray*}
Note that $v_{kl}\in H^{j}(\R_{+},r^{N-1}\D r)$ holds for all $k,l\ge0$
since $v\in H^{j}(\R^{N})$.

Now we start the proof of Theorem \ref{thm: Nondegeneracy}. We first
prove that $u$ is nondegenerate in $H_{\rad}^{1}(\R^{3})$ (in the
sense of the following proposition), which is the key ingredient of
the proof of Theorem \ref{thm: Nondegeneracy}.

\begin{proposition} \label{prop: nonexistence of radial solutions}
Let $\L_{+}$ be defined as in Definition \ref{def: Nondegeneracy}
and let $\var\in H_{\rad}^{1}(\R^{3})$ be such that $\L_{+}\var=0$.
Then $\var\equiv0$ in $\R^{3}$. \end{proposition}
\begin{proof}
Let $\var\in H_{\rad}^{1}(\R^{3})$ be such that $\L_{+}\var=0$.
By virtue of the notations introduced above, we can rewrite the equation
$\L_{+}\var=0$ as below:
\[
{\cal A}_{u}\var=2b\left(\int_{\R^{3}}\na u\cdot\na\var\D x\right)\De u.
\]
We have to prove that $\var\equiv0$. This is sufficient to show that
\begin{equation}
\int_{\R^{3}}\na u\cdot\na\var\D x=0,\label{eq: Key observation}
\end{equation}
since then $\var\in\Ker{\cal A}_{u}\cap L_{\rad}^{2}(\R^{3})$, which
implies that $\var\equiv0$ by Proposition \ref{prop: property of Au}.

To deduce (\ref{eq: Key observation}), we proceed as follows. Since
$u$ is radial and ${\cal A}_{u}$ is one-to-one on $L_{\rad}^{2}(\R^{3})$
by Proposition \ref{prop: property of Au}, $\var$ satisfies the
equivalent equation
\[
\var=2b\left(\int_{\R^{3}}\na u\cdot\na\var\D x\right){\cal A}_{u}^{-1}(\De u),
\]
where ${\cal A}_{u}^{-1}$ is the inverse of ${\cal A}_{u}$ restricted
on $L_{\rad}^{2}(\R^{3})$.

Next we compute ${\cal A}_{u}^{-1}(\De u)$. By Eq. (\ref{eq: auxiliary equation}),
$\De u=(u-u^{p})/c$. Hence ${\cal A}_{u}^{-1}(\De u)=\left({\cal A}_{u}^{-1}(u)-{\cal A}_{u}^{-1}(u^{p})\right)/c$.
Applying Proposition \ref{prop: property of Au} (3), we deduce that
\[
{\cal A}_{u}^{-1}(\De u)=\frac{1}{c}\left(-\frac{S}{2}+\frac{u}{p-1}\right)=-\frac{1}{2c}x\cdot\na u,
\]
where $S$ is defined as in Proposition \ref{prop: property of Au}.
Therefore, we obtain
\[
\var=-\frac{b}{c}\left(\int_{\R^{3}}\na u\cdot\na\var\D x\right)x\cdot\na u=\left(\int_{\R^{3}}\na u\cdot\na\var\D x\right)\psi,
\]
with $\psi=-\frac{b}{c}x\cdot\na u$

Now we can deduce (\ref{eq: Key observation}) from the above
formula. Taking gradient on both sides gives
\[
\na\var=\left(\int_{\R^{3}}\na u\cdot\na\var\D x\right)\na\psi.
\]
Multiply $\na u$ on both sides and integrate. We achieve
\[
\int_{\R^{3}}\na u\cdot\na\var\D x=\left(\int_{\R^{3}}\na u\cdot\na\var\D x\right)\int_{\R^{3}}\na u\cdot\na\psi\D x.
\]
A direct computation yields that
\[
\int_{\R^{3}}\na u\cdot\na\psi\D x=\frac{b}{2c}\int_{\R^{3}}|\na u|^{2}\D x=\frac{c-a}{2c}<\frac{1}{2}.
\]
Hence we easily deduce that $\int_{\R^{3}}\na u\cdot\na\var\D x=0$,
that is, (\ref{eq: Key observation}) holds. The proof of Proposition
\ref{prop: nonexistence of radial solutions} is complete.
\end{proof}
With the help of Proposition \ref{prop: nonexistence of radial solutions}, we can now finish the proof of Theorem \ref{thm: Nondegeneracy}. The procedure is standard, see  e.g. Ambrosetti and Malchiodi
\cite[Section 4.2]{Ambrosetti-Malchiodi-Book}. For the readers' convenience, we  give a detailed proof.

\begin{proof}[Proof of Theorem \ref{thm: Nondegeneracy}] Let $\var\in
H^{1}(\R^{3})$ be such that $\L_{+}\var=0$. We have to prove that
$\var$ is a linear combination of $\pa_{x_{i}}u$, $i=1,2,3$. The idea is to turn the problem $\L_{+}\var=0$ into a system of ordinary differential equations by making
use of the spherical harmonics to decompose $\var$ into
\[
\var=\sum_{k=0}^{\wq}\sum_{l=1}^{M_{k}-M_{k-2}}\var_{kl}(r)Y_{kl}(\Om)
\]
with $r=|x|$ and $\Om=x/|x|$, where
\begin{eqnarray}
\var_{kl}(r)=\int_{\S^{2}}\var(r\Om)Y_{kl}(\Om)\D\Om &  & \forall\, k\ge0.\label{eq: kl-th component}
\end{eqnarray}
Note that $\var_{kl}\in H^{1}(\R_{+},r^{2}\D r)$ holds for all
$k,l\ge0$ since $\var\in H^{1}(\R^{3})$.

Combining the fact that  $\int_{\S^2}Y_{kl}\D\si=0$ hold for
all $k,l\ge1$,  together with the fact that $u$
is radial, we deduce
\[
\int_{\R^{3}}\na u\cdot\na\var\D x=\int_{\R^{3}}\left(-\De
u\right)\var\D x=\int_{\R^{3}}\na u\cdot\na\var_{0}\D x,
\]
where $\var_0(x)=\var_{0,1}(|x|)$ for $x\in\R^3$.

 Hence, the problem
$\L_{+}\var=0$ is equivalent to the following system of ordinary
differential equations:

For $k=0$, we have
\begin{equation}
\L_{+}\var_{0}=0.\label{eq: first case}
\end{equation}

For $k=1$, we have
\begin{equation}
A_1(\var_{1l})\equiv\left(-c\De_{r}+\frac{\la_{1}}{r^{2}}\right)\var_{1l}+\var_{1l}-pu^{p-1}\var_{1l}=0\label{eq:
second case}
\end{equation}
for $l=1,2,3$. Here $\De_r=\pa_{rr}+\frac 2r \pa_r$. We also used
the fact that $u$ and $\De u$ are radial functions.

For $k\ge2$, we have that
\begin{equation}
A_k(\var_{kl})\equiv\left(-c\De_{r}+\frac{\la_{k}}{r^{2}}\right)\var_{kl}+\var_{kl}-pu^{p-1}\var_{kl}=0.\label{eq:
third case}
\end{equation}

To solve  Eq. (\ref{eq: first case}), we apply Proposition \ref{prop: nonexistence
of radial solutions} to conclude that $\var_{0}\equiv0$.

To solve  Eq. (\ref{eq: second case}),  note
 that $u^{\prime}$ is a solution of Eq. (\ref{eq: second case})  and
$u^{\prime}\in H^{1}(\R_{+},r^{2}\D r)$. Since  Eq. (\ref{eq: second case}) is a  second order linear ordinary differential equation, we assume that it has another solution $v(r)=h(r)u^{\prime}(r)$ for some $h$. It is easy to find that $h$ satisfies
$$h^{\prime\prime}u^{\prime}+\frac {2}{r}h^{\prime}u^{\prime}+2h(u^{\prime})^{\prime}=0.$$
If $h$ is not identically a constant,  we derive that
$$-\frac {h^{\prime\prime}}{h^{\prime}}=2\frac {u^{\prime\prime}}{u^{\prime}}+\frac 2r,$$
which implies that
\begin{eqnarray*}
 h^{\prime}(r)\sim r^{-2}(u^{\prime})^2&&\text{as}\,\,r\to \infty.
\end{eqnarray*}
Recall that   $Q=Q(|x|)$, $x\in \R^3$, is the unique positive radial solution of Eq. (\ref{eq: Kwong}).
It is well known \cite{GNN1981} that $\lim_{r\to \infty}re^rQ^{\prime}(r)= -C$ holds  for some constant $C>0$. Hence, by the proof of Theorem \ref{thm: Uniqueness}, we know that $\lim_{r\to \infty}re^{r/\sqrt{c}}u^{\prime}(r)=-C_1$ for some $C_1>0$. Combining this fact with the above estimates gives
$$|h(r)u^{\prime}(r)|\ge Cr^{-1}e^{r/\sqrt{c}}$$
as $r\to \infty$. Thus $hu^{\prime}$ does not belong to $H^1(\R_+,r^2\D r)$ unless $h$ is a constant. This shows that the family of solutions of  Eq. (\ref{eq: second case}) in  $H^1(\R_+,r^2\D r)$ is given by $hu^{\prime}$, for some constant $h$.
In particular, we conclude that
$\var_{1l}=d_{l}u^{\prime}$ hold for some constant $d_{l}$, for all $1\le l\le 3$.

For the last Eq.  (\ref{eq:
third case}),   we show that it has only a trivial solution.
Indeed, for  $k\ge2$, we have
\[A_k=A_1+\frac {\de_k}{r^2},\]
where $\de_k=\la_k-\la_1$. Since $\la_{k}>\la_{1}$, we find that $\de_k >0$.  Notice that $u^{\prime}$ is an
eigenfunction of $A_1$ corresponding to the eigenvalue 0, and that
$u^{\prime}$ is of constant sign. By virtue of orthogonality, we
can easily infer that $0$ is the smallest eigenvalue of $A_1$.
That is, $A_1$ is a nonnegative operator. Therefore, $\de_k>0$
implies that $A_k$ is a  positive operator for all $k\ge 2$. That is, $\langle A_k\psi,\psi \rangle \ge 0$ for all $\psi \in H^1(\R_+,r^2 \D r)$, and the equality attains if and only if $\psi =0$. As a result, we easily prove that if $\var_{kl}$ is a solution of Eq. (\ref{eq: third case}), then
$\var_{kl}\equiv0$ holds for all $k\ge 2$.

In summary, we obtain
\[
\var=\sum_{l=1}^{3}d_{l}u^{\prime}(r)Y_{1l}=\sum_{l=1}^{3}d_{l}\pa_{x_{l}}u.
\]
The proof of Theorem \ref{thm: Nondegeneracy} is complete. \end{proof}



\begin{thebibliography}{10}


\bibitem{Shinji-Masataka-Tatsuya-2016} \textsc{S. Adachi, M. Shibata and T. Watanabe},
\emph{A note on the uniqueness and the non-degeneracy of positive
radial solutions for semilinear elliptic problems and its application.
}Preprint at arXiv: arXiv:1602.07086 {[}math.AP{]}.

\bibitem{Ambrosetti-Malchiodi-Book} \textsc{A. Ambrosetti and A. Malchiodi,}
\emph{Perturbation Methods and Semilinear Elliptic Problems on }$\R^{N}$.
Birkh\"{a}user Verlag, 2006.

\bibitem{DAncona-Spagnolo-1992} \textsc{P. D'Ancona and S. Spagnolo,}
\emph{Global solvability for the degenerate Kirchhoff equation with
real analytic data}. Invent. Math. \textbf{108} (1992), 247-262.

\bibitem{Arosio-Panizzi-1996} \textsc{A. Arosio and S. Panizzi,}
\emph{On the well-posedness of the Kirchhoff string.} Trans. Amer.
Math. Soc. \textbf{348} (1996), 305-330.


\bibitem{Berestycki-Lions-1983-1} \textsc{H. Berestycki and P.-L. Lions},
\emph{Nonlinear scalar field equations. I. Existence of a ground state.}
Arch. Rational Mech. Anal. \textbf{82} (1983), no. 4, 313-345.


\bibitem{Bernstein-1940} \textsc{S. Bernstein,} \emph{Sur une classe
d'\'equations fonctionelles aux d\'eriv\'ees partielles.} Bull.
Acad. Sci. URSS. S\emph{\'e}r. \textbf{4} (1940), 17-26.


\bibitem{Chang et al-2007}  \textsc{S.-M. Chang, S. Gustafson, K. Nakanishi and T.-P. Tsai,}
\emph{Spectra of linearized operators for NLS solitary waves.} SIAM
J. Math. Anal. \textbf{39} (2007/08), no. 4, 1070-1111.

\bibitem{Cingolani-Lazzo-1997} \textsc{S. Cingolani and N. Lazzo,}\emph{
Multiple semiclassical standing waves for a class of nonlinear Schr\"odinger
equations.} Topol. Methods Nonlinear Anal. \textbf{10} (1997), 1-13.



\bibitem{Deng-Peng-Shuai-2015} \textsc{Y. Deng, S. Peng and W. Shuai,}
\emph{Existence and asymptotic behavior of nodal solutions for the
Kirchhoff-type problems in $\R^{3}$.} J. Funct. Anal. \textbf{269}
(2015), no. 11, 3500-3527.

\bibitem{Fall-Valdinoci-2014} \textsc{M. M. Fall and E. Valdinoci,}
\emph{Uniqueness and nondegeneracy of positive solutions of $(-\De)^{s}u+u=u^{p}$
in $\R^{N}$ when s is close to 1.} Comm. Math. Phys. \textbf{329}
(2014), no. 1, 383-404.

\bibitem{Figueiredo et al-2014} \textsc{G.M. Figueiredo, N. Ikoma, N. and J. R.  Santos J\'unior,}\emph{
Existence and concentration result for the Kirchhoff type equations
with general nonlinearities.} Arch. Rational Mech. Anal. \textbf{213}
(2014), 931-979.

\bibitem{Frank-Lenzmann-2013}  \textsc{R.L. Frank and E. Lenzmann,}
\emph{Uniqueness of non-linear ground states for fractional Laplacians
in $\R$.} Acta Math. \textbf{210} (2013), no. 2, 261-318.

\bibitem{Frank-Lenzman-Silvestre-2016} \textsc{R.L. Frank, E.  Lenzmann and L. Silvestre,}
\emph{Uniqueness of radial solutions for the fractional Laplacian.}
Commun. Pur. Appl. Math. \textbf{69} (2016), 1671-1726.

\bibitem{GNN1981}\textsc{B. Gidas, W.M. Ni and L. Nirenberg},  \emph{Symmetry of positive solutions of nonlinear elliptic equations in $\R^{n}$.} Mathematical analysis and applications, Part A, pp. 369-402, Adv. in Math. Suppl. Stud., 7a, Academic Press, New York-London, 1981.

\bibitem{Guo-2015} \textsc{Z. Guo}, \emph{Ground states for Kirchhoff
equations without compact condition.} J. Differential Equations \textbf{259}
(2015), no. 7, 2884-2902.

\bibitem{He-2016} \textsc{Y. He}, \emph{Concentrating bounded states for a class of singularly perturbed Kirchhoff type equations with a general nonlinearity.} J. Differential Equations \textbf{261}
(2016),  6178-6220.

\bibitem{He-Li-2015} \textsc{Y. He and G. Li,} \emph{Standing waves
for a class of Kirchhoff type problems in $\R^{3}$ involving critical
Sobolev exponents.} Calc. Var. Partial Differential Equations \textbf{54}
(2015), no. 3, 3067-3106.

\bibitem{He-Li-Peng-2014} \textsc{Y. He,  G. Li and S. Peng}, \emph{Concentrating
bound states for Kirchhoff type problems in $\R^{3}$ involving critical
Sobolev exponents.} Adv. Nonlinear Stud. \textbf{14} (2014), no. 2,
483-510.

\bibitem{He-Zou-2009} \textsc{X. He and  W. Zou,} \emph{Infinitely
many positive solutions for Kirchhoff-type problems}. Nonlinear Anal.
\textbf{70 }(2009), 1407-1414.

\bibitem{He-Zou-2010} \textsc{X. He and  W. Zou,} \emph{Multiplicity
of solutions for a class of Kirchhoff type problems}. Acta Math. Appl.
Sin. \textbf{26} (2010) 387-394.

\bibitem{He-Zou-2012} \textsc{X. He and  W. Zou,} \emph{Existence
and concentration behavior of positive solutions for a Kirchhoff equation
in ${\cal R}^{3}$.} J. Differential Equations \textbf{252} (2012),
1813-1834.


\bibitem{Kirchhoff-1883} \textsc{G. Kirchhoff,} Mechanik, Teubner,
Leipzig, 1883.

\bibitem{Kwong-1989} \textsc{M. K. Kwong,} \emph{Uniqueness of positive
solutions of $\De u-u+u^{p}=0$ in $\mathbf{R}^{n}$.} Arch. Rational
Mech. Anal. \textbf{105} (1989), no. 3, 243-266.


\bibitem{Lenzmann-2009} \textsc{E. Lenzmann,} \emph{Uniqueness of
ground states for pseudorelativistic Hartree equations.} Anal. PDE
\textbf{2} (2009), no. 1, 1-27.

\bibitem{Li-Li-Shi-2012} \textsc{Y. Li, F. Li and J. Shi,} \emph{Existence
of a positive solution to Kirchhoff type problems without compactness
conditions.} J. Differential Equations \textbf{253} (2012), 2285-2294.

\bibitem{Li-Ye-2014} \textsc{G. Li and H. Ye,} \emph{Existence of
positive ground state solutions for the nonlinear Kirchhoff type equations
in $\R^{3}$.} J. Differential Equations \textbf{257} (2014), no.
2, 566-600.

\bibitem{Lions-1978} \textsc{J.L. Lions,} \emph{On some questions
in boundary value problems of mathematical physics}. Contemporary
Development in Continuum Mechanics and Partial Differential Equations,
in: North-Holland Math. Stud., vol. 30, North-Holland, Amsterdam,
New York, 1978, pp. 284-346.

\bibitem{Lu-2015-arXiv} \textsc{S.-S. Lu}, \emph{An autonomous Kirchhoff-type
equation with general nonlinearity in $\R^{N}$.} Preprint at arXiv:
arXiv:1510.07231v2 {[}math.AP{]}.


\bibitem{Perera-Zhang-2006} \textsc{K. Perera and Z. Zhang,} \emph{Nontrivial
solutions of Kirchhoff-type problems via the Yang index}. J. Differential
Equations \textbf{221} (2006), 246-255.

\bibitem{Pohozaev-1975} \textsc{S.I. Pohozaev}, \emph{A certain
class of quasilinear hyperbolic equations.} Mat. Sb. (N.S.) \textbf{96}
(138) (1975), 152-166, 168 (in Russian).

\bibitem{Selvitella-2015}  \textsc{A. Selvitella}, \emph{Nondegeneracy
of the ground state for quasilinear Schr\"odinger equations. } Calc.
Var. Partial Differential Equations\textbf{ 53 }(2015), 349-364.


\bibitem{Shuai-2015} \textsc{W. Shuai}, \emph{Sign-changing solutions
for a class of Kirchhoff-type problem in bounded domains.} J. Differential
Equations \textbf{259} (2015), no. 4, 1256-1274.

\bibitem{Xiang-2016} \textsc{C.-L. Xiang}, \emph{Remarks on Nondegeneracy
of Ground States for Quasilinear Schr\"odinger Equations. } Discrete
Contin. Dyn. Syst. \textbf{36} (2016), no. 10, 5789-5800.

\bibitem{Xiang-CVPDE} \textsc{C.-L. Xiang}, \emph{Uniqueness and Nondegeneracy of Ground States for Choquard Equations in Three Dimensions.} To apppear in Calc.
Var. Partial Differential Equations.

\bibitem{Ye-2015} \textsc{H. Ye}, \emph{Positive high energy solution
for Kirchhoff equation in $\R^{3}$ with superlinear nonlinearities
via Nehari-Poho\u{z}aev manifold.} Discrete Contin. Dyn. Syst. \textbf{35}
(2015), no. 8, 3857-3877.

\bibitem{Zhang-Perera-2006} \textsc{Z. Zhang and K. Perera,} \emph{Sign
changing solutions of Kirchhoff type problems via invariant sets of
descent flow.} J. Math. Anal. Appl. \textbf{317} (2006), no. 2, 456-463.

\bibitem{Ziemer} \textsc{W.P. Ziemer}, \emph{Weakly differentiable
functions.} Graduate Texts in Mathematics, 120. Springer-Verlag, New
York, 1989. \end{thebibliography}
\end{document}